\documentclass[oneside, 12pt]{amsart}

\usepackage{amscd,amssymb,enumerate, amsmath,graphicx}
\usepackage{mathrsfs}
\usepackage[all,cmtip]{xy}

\newcommand\mylabel[1]{\label{#1}}

\newtheorem{theorem}{Theorem}[section]
\newtheorem*{maintheorem}{Theorem}
\newtheorem{lemma}[theorem]{Lemma}
\newtheorem{proposition}[theorem]{Proposition}

\theoremstyle{definition}

\theoremstyle{remark}

\DeclareFontFamily{U}{wncy}{}
\DeclareFontShape{U}{wncy}{m}{n}{<->wncyr10}{}
\DeclareSymbolFont{mcy}{U}{wncy}{m}{n}
\DeclareMathSymbol{\Sh}{\mathord}{mcy}{"58}
\newcommand{\ZZ}	{\mathbb{Z}}
\newcommand{\QQ}	{\mathbb{Q}}
\newcommand{\RR}	{\mathbb{R}}

\newcommand{\FF}	{\mathbb{F}}

\newcommand{\PP}	{\mathbb{P}}

\newcommand{\ideala}	{\mathfrak{a}}
\newcommand{\idealb}	{\mathfrak{b}}

\newcommand  {\shE}     {\mathscr{E}}
\newcommand  {\shF}     {\mathscr{F}}

\newcommand  {\shI}     {\mathscr{I}}

\newcommand  {\shM}     {\mathscr{M}}

\newcommand  {\shN}     {\mathscr{N}}
\newcommand  {\shL}     {\mathscr{L}}

\newcommand  {\Ass}     {\operatorname{Ass}}

\newcommand  {\cH}      {\check{H}}

\newcommand  {\Coh}   	{\operatorname{Coh}}

\newcommand  {\Et}      {{\text{\rm Et}}}

\newcommand  {\Frac}    {\operatorname{Frac}}

\newcommand  {\id}      {{\operatorname{id}}}

\newcommand  {\length}	{\operatorname{length}}

\newcommand  {\dirlim}  {\varinjlim}
\newcommand  {\invlim}  {\varprojlim}

\newcommand  {\lra}     {\longrightarrow}

\newcommand  {\maxid}   {\mathfrak{m}}

\newcommand  {\mult}    {\operatorname{mult}}

\newcommand  {\Nil}     {\operatorname{Nil}}
\newcommand  {\Nis}  	{\text{\rm Nis}}

\newcommand  {\primid}  {\mathfrak{p}}
\renewcommand{\O}       {\mathscr{O}}

\newcommand  {\pr}      {\operatorname{pr}}
\newcommand  {\Proj}    {\operatorname{Proj}}

\newcommand  {\quadand} {\quad\text{and}\quad}

\newcommand  {\ra}      {\rightarrow}

\newcommand  {\red}     {{\operatorname{red}}}

\newcommand  {\sep}     {{\operatorname{sep}}}
\newcommand  {\Set}     {{\text{\rm Set}}}

\newcommand  {\Spec}    {\operatorname{Spec}}

\newcommand  {\Sym}     {\operatorname{Sym}}

\newcommand  {\TSC}      {\operatorname{TSC}}

\newcommand  {\vol}  	{\operatorname{vol}}

\def\mydate{\number\day\space\ifcase\month \or January\or February\or March\or 
April\or May\or June\or July\or
August\or September\or October\or November\or December\fi \space\number\year}

\newcommand{\tX}{\tilde{X}}
\newcommand{\tu}{\tilde{u}}
\newcommand{\tv}{\tilde{v}}
\newcommand{\tA}{\tilde{A}}
\newcommand{\tB}{\tilde{B}}
\newcommand{\tC}{\tilde{C}}
\newcommand{\tP}{\tilde{P}}

\newcommand{\tshL}{\tilde{\shL}}
\renewcommand{\emptyset}{\varnothing}
\newcommand{\catC}{\mathcal{C}}

\DeclareFontFamily{U}{wncy}{}
\DeclareFontShape{U}{wncy}{m}{n}{<->wncyr10}{}
\DeclareSymbolFont{mcy}{U}{wncy}{m}{n}
\DeclareMathSymbol{\Sh}{\mathord}{mcy}{"58}
\begin{document}

\title[Total separable closure]
      {Total separable closure and contractions}

\author[Stefan Schr\"oer]{Stefan Schr\"oer}
\address{Mathematisches Institut, Heinrich-Heine-Universit\"at,
40204 D\"usseldorf, Germany}
\curraddr{}
\email{schroeer@math.uni-duesseldorf.de}


\dedicatory{22 August 2017}

\begin{abstract}
We show that on  integral normal separated schemes whose function field
is separably closed, for each pair of points the intersection of
the resulting local schemes is local. This extends a result of Artin
from rings to   schemes. The argument relies
on the existence of certain modifications in inverse limits.
As an application, we show that \v{C}ech cohomology coincides with sheaf
cohomology for the Nisnevich topology. Along the way, we generalize
the characterization of contractible curves on surfaces by 
negative-definiteness of the  intersection matrix to higher dimensions,
using bigness of invertible sheaves on non-reduced schemes.
\end{abstract}

\maketitle
\tableofcontents

%===========================================================
\section*{Introduction}
\mylabel{Introduction}

This paper deals with the  Zariski topology for a  class of  schemes that are in general highly non-noetherian,
yet arise from noetherian schemes in a canonical way:
We say  that an integral   scheme $X$ is   \emph{totally separably closed} if it
is normal and the function field $F=\O_{X,\eta}=\kappa(\eta)$ is separably closed.
As abbreviation  one also  says that $X$ is a \emph{TSC scheme}.
Each integral scheme $X_0$ has a total separable closure $X=\TSC(X_0)$, defined
as the integral closure with respect to a chosen separable closure  $F=F_0^\sep$ of the function field $F_0$.
Such schemes $X$ are \emph{everywhere strictly local}. In other words,
all local rings $\O_{X,x}$, $x\in X$ are strictly local rings, that is, henselian with separably closed residue field. 
One may regard them as analogues of \emph{Pr\"ufer schemes}, where
all local rings are valuation rings.

TSC schemes have some relevance with respect to the \'etale topology. Indeed, M.\ Artin \cite{Artin 1971}
used them to prove that    \v{C}ech cohomology equals sheaf cohomology
for the \'etale topology over affine schemes $X=\Spec(R)$. This result
immediately extends to schemes with the \emph{AF property}, which means that
any finite subset admits an affine open neighborhood.
Note, however, that by  \cite{Benoist 2013}, Corollary 2, the AF property is equivalent to quasiprojectivity for
normal schemes that are separated and of finite type over a ground field.
Actually, Artin used algebraic closure rather that separable closure, but this  makes no
difference for the underlying topological spaces.
See Huneke's overview \cite{Huneke 2011} for the role of absolute integral closure
in commutative algebra.

One crucial step in Artin's arguments is to show that affine integral  TSC schemes $X$
have the following  surprising property, which is of purely topological nature: For any pair of points $u,v\in X$ the
intersection of local schemes
$$
\Spec(\O_{X,u})\cap\Spec(\O_{X,v})\subset X
$$
remains  a local scheme. If we endow  the underlying set $X$ with the order relation 
$x\leq y \Leftrightarrow x\in\overline{\{y\}}$, the above property means that the supremum
$\sup(u,v)$ exists for all pairs of points $u,v\in X$. 
This strange property almost never holds on  noetherian schemes $X_0$,
and intuitively means that    in    inverse limits $X=\invlim X_\lambda$,    
common generizations of $u_\lambda,v_\lambda\in X_\lambda$ 
are totally  ``ripped apart''.
In some sense, this is a topological incarnation of the result of Schmidt that 
a field with two different henselian valuations is separably closed
(\cite{Schmidt 1933}, Satz 3. See   \cite{Engler and Prestel 2005}, Theorem 4.4.1 for a modern account).
The main goal of this paper is to establish Artin's result in full generality:

\begin{maintheorem}{\bf (See Theorem \ref{local})}
For any separated integral TSC scheme $X$, the   intersections 
$\Spec(\O_{X,u})\cap\Spec(\O_{X,v})\subset X$
are local for all points $u,v\in X$.
\end{maintheorem}

In \cite{Schroeer 2017}, I already obtained this   for   total separable closures
of schemes $X_0$ that are separated and of finite type over a ground field $k$.
The arguments rely on modifications and contractions in inverse limits $X=\invlim X_\lambda$,
and do not apply in mixed characteristics. Here we modify our approach,
and reduce the problem to proper schemes over  excellent Dedekind domains.
We then use different modifications $X'$ and contractions $\tX$ in inverse limits 
so that   Artin's result applies to the  TSC scheme $\tX$, which is constructed to have the AF property.
This is enough to conclude for the original TSC scheme $X$.

To carry this out, we have to analyze the existence of suitable \emph{modifications} and \emph{contractions}.
On     algebraic surfaces $X$, a curve  $E=E_1+\ldots+E_r$ is contractible to points 
if and only if the intersection matrix $\Phi=(E_i\cdot E_j)$ is negative-definite.
This observation goes back to Mumford, Artin and Deligne, in various forms of generality.
Note that in general the   contractions $r:X\ra Y$ yield  algebraic spaces rather than  schemes.
The following generalization to higher dimensions seems to be of independent interest:

\begin{maintheorem}{\bf (See Theorem \ref{big})}
Let $X$ be  a normal scheme that is proper over an excellent Dedekind domain $R$, and $E=E_1+\ldots+E_r$ be a Weil divisor contained in
a closed fiber for the structure morphism $X\ra\Spec(R)$. If $E$ is  is
contractible to points, then for each effective Cartier divisor $D=\sum m_iE_i$, the   invertible sheaf
$\O_D(-D)$ on $D$ is big.
\end{maintheorem}

Here bigness for an invertible sheaf $\shL$ on some proper algebraic scheme $Z$, which is not necessarily reduced or irreducible,
is defined in terms of the \emph{Iitaka dimension}, which itself is given, up to  a shift, by the 
Krull dimension of the ring $R(Z,\shL)=\bigoplus_{n\geq 0} H^0(Z,\shL^{\otimes n})$.
This   generalization from integral to arbitrary schemes was  analyzed by Cutkosky \cite{Cutkosky 2014},
and his results on the multiplicity or volume $\mult(\ideala_\bullet)=\vol(\ideala_\bullet)$
for \emph{graded families of ideals}, together with Huneke's version \cite{Huneke 1992} of the Brian\c{c}on--Skoda Theorem,   play a crucial role for the above.

As explained in  \cite{Schroeer 2017}, our main result on TSC schemes has  immediate
consequences for the \emph{Nisnevich topology} of completely decomposed \'etale maps
\cite{Nisnevich 1989}. This is a variant of the \'etale topology, where the local rings are 
  henselian local rings rather than   strictly local rings. We get:

\begin{maintheorem}{\bf (See Theorem \ref{cech})}
For each quasicompact separated scheme $X$ and every abelian Nisnevich sheaf $F$, the canonical maps
$$
\check{H}^p_\Nis(X,F)\lra H^p_\Nis(X,F) 
$$
from \v{C}ech cohomology to sheaf cohomology are bijective in every degree $p\geq 0$.
\end{maintheorem}

The paper is organized as follows:
Section 1 contains our results on contractions for proper schemes over excellent Dedekind domains.
The main result about TSC schemes is given in Section 2.
The final Section 3 gives the application to   Nisnevich cohomology.

%===========================================================
\section{Contractions over Dedekind domains}
\mylabel{Contractions}

Let $S=\Spec(R)$ be the spectrum of a Dedekind domain $R$, and $X$ be an  proper $S$-scheme.
We write $f:X\ra S$ for the structure morphism. For simplicity, we assume
that the scheme $X$ is integral and that the ring $R$ is excellent. Note that we do not assume flatness; in particular, the
structure morphism may factor over some closed point $\sigma\in S$.

A closed subscheme $E\subset X$ is said to be \emph{contractible to points}
if there is a commutative diagram
$$
\xymatrix{
X\ar[rr]^r\ar[dr]_f	&	& Y\ar[dl]^g\\
		& S	&\\		
}
$$
where $Y$ is an algebraic space, the structure morphism $g:Y\ra S$ is proper,
and $r:X\ra Y$ is a morphism with $\O_Y=r_*(\O_X)$ that is an open embedding
on $X\smallsetminus E$ such that the image $Z=r(E)$ consists of finitely many closed points.
Their images in $S$ are closed as well, and it follows that the connected components
of $E$ are contained in   closed fibers for the structure morphism $f:X\ra S$.
The morphism $r:X\ra Y$ is unique up to unique isomorphism,
and depends only on the underlying closed set for $E\subset X$, which follows
from  \cite{EGA II}, Lemma 8.11.1.

Algebraic spaces can be glued along open   subsets in the same way as ringed spaces
(a consequence  from \cite{Olsson 2016}, Proposition 5.2.5).
In particular,  the closed subset $E\subset X$ is contractible to points if and only if each
connected component is \emph{contractible to a single point}. 
Using Stein factorization, one easily sees  the following permanence property:

\begin{proposition}
\mylabel{pullback}
Assume  $X'$ is another proper $S$-scheme that is integral,
and let $X'\ra X$ be a     morphism. 
Suppose that  a closed subset $E\subset X$ is contractible to points, and that the morphism $X'\ra X$
is finite over $X\smallsetminus E$. Then the preimage 
 $E'=E\times_XX'$ is contractible to points. Indeed, if $X\ra Y$ is the contraction of $E\subset Y$,
then the Stein factorization  $Y'$ for the composition $X'\ra Y$ is the contraction
of $E'\subset X'$.
\end{proposition}

For general closed subsets $E\subset X$ it is often  difficult
to verify contractibility. However, by applying the previous result to the blowing-up $X'\ra X$ with center $E$
one reduces to the case of effective Cartier divisors.
Then more can be said:

\begin{proposition}
\mylabel{contractible}
Suppose  $E\subset X$ is an effective Cartier divisor contained in some closed fiber $X_\sigma=f^{-1}(\sigma)$.
Let $\shL=\O_X(-E)$. If the restriction 
$\shL|E$ is ample  then the closed   subset $E\subset X$ is contractible to   points.
\end{proposition}

\proof
This  immediately follows from  Corollary 6.10 in  Artin's work  \cite{Artin 1970} on algebraic stacks:
It suffices to treat the case that $E$ is connected, and 
we need to check two conditions.
The first condition is straightforward: for every coherent sheaf $\shF$ on $E$ the cohomology group $H^1(E,\shF\otimes\shL^{\otimes n})$ vanishes
for all $n\gg 0$, because $\shL|E$ is ample.
The second conditions is somewhat more intricate:
Since $E$ is proper, connected and contained in a closed fiber $X_\sigma=f^{-1}(\sigma)$, the rings 
$R_n=H^0(X,\O_{nE})$ are finite local artinian $R$-algebras. Write $k_n=R_n/\maxid_{R_n}$ for 
their residue fields. The inclusions $E\subset 2E\subset\ldots$ induce an an inverse system
$k_1\supset k_2\supset\ldots$ of fields, all of which contain the residue field $\kappa=\kappa(\sigma)$ and have
finite degree.
 Let $k=\bigcap k_n$ be their intersection,
and choose an index $n_0$ so that the inclusions $k_{n+1}\subset k_n$ are equalities
for all $n\geq n_0$.
Consider the resulting morphism $E\ra\Spec(k)$ and, for each $n\geq 0$,   the cartesian diagram
$$
\begin{CD}
R_n\times_{k_1} k	@>>>	k\\
@VVV				@VVV\\
R_n			@>>>	k_1.
\end{CD}
$$
Artin's second condition stipulates that the upper vertical arrows must be surjective.
To see this,  choose some index $m\geq \max(n_0,n)$. Then $k_m=k$ and the residue class map $R_m\ra k_m=k$ factors over
the fiber product $R_n\times_{k_1}k$, so the projection in question is surjective.
\qed

\medskip
The  converse   does not hold: For example, if $X$ is regular of dimension $d=2$,
and $E=E_1+E_2$ is a   curve with two irreducible components, having  intersection matrix $\Phi=(\begin{smallmatrix}-2&1\\1&-2 \end{smallmatrix})$.
The latter is negative-definite, so the curve $E$ is contractible.
The linear combination $D=3E_1+E_2$ is contractible as well, yet $\O_D(-D)$ is not ample, because
$(D\cdot E_2)=-1$.

Nevertheless, it is   natural to ask for some   form of   converse.
Indeed, we shall establish such a result based on the notion of  bigness  rather then ampleness.
Let us  recall the relevant definitions:
Suppose  $Z$ is a proper scheme over some ground field $k$.
Given an invertible sheaf $\shL$ on $Z$ we get
a graded ring $R(Z,\shL)=\bigoplus_{n\geq 0} H^0(Z,\shL^{\otimes n})$,
which is is not necessarily of finite type or noetherian. Let $d=\dim(R)$ be its Krull dimension. The \emph{Iitaka dimension}
or \emph{Kodaira--Iitaka dimension}  is defined as
$$
\kappa(\shL)=\begin{cases}
d-1 		& \text{if $d\geq 1$;}\\
-\infty		& \text{else.}
\end{cases}
$$
Note that  it will be crucial to allow reducible and non-reduced schemes $Z$ for what we have in mind.
For integral normal schemes $Z$, the Iitaka dimension is a classical notion from 
birational geometry: If some $\shL^{\otimes n_0}$ with $n_0\geq 1$ has a non-zero global section, 
the number $\kappa(\shL)$   can also be seen as the maximal dimension of the images for 
the rational maps $X\dashrightarrow \PP^m$ defined by $\shL^{\otimes n}$, where $m=h^0(\shL^{\otimes n})-1$ and $n\geq 1$ runs over
the positive multiplies of $n_0$. We refer to the monograph  of Lazarsfeld (\cite{Lazarsfeld 2004}, Section 2.2)  for more details.
Iitaka dimension  was only recently extended to arbitrary proper schemes, by the work
of Cutkosky on asymptotics of 
ideals and linear series. In fact, in \cite{Cutkosky 2014}, Section 7 he defined it in the more general context
of \emph{graded linear series}, which can be seen as graded  subrings  $L=\bigoplus_{n\geq 0} L_n$ inside $R(Z,\shL)$.

According to \cite{Cutkosky 2014}, Lemma 7.1 we have $\kappa(\shL)\leq \dim(Z)$   for
arbitrary proper schemes $Z$.
In case of equality $\kappa(\shL)=\dim(Z)$ one says that the invertible sheaf $\shL$ is \emph{big}.
We need the following observation on invertible sheaves that are not big:

\begin{lemma}
\mylabel{devissage}
Set $d=\dim(Z)$. If $\shL$ is not big, then  for each coherent
sheaf $\shF$ on $Z$ there is a real constant $\beta\geq 0$ so that
$h^0(\shF\otimes\shL^{\otimes n})\leq \beta n^{d-1}$ for all integers $n\geq 0$.
\end{lemma}

\proof
This is a devissage argument similar to \cite{EGA IIIa}, Theorem 3.1.2.
Let $\Coh(Z)$ be the abelian category of all coherent sheaves $\shF$ on $Z$, and  $\catC\subset\Coh(Z)$ be the   subcategory of all sheaves
for which the assertion holds. If  $0\ra\shF'\ra\shF\ra\shF''\ra 0$ is a short exact sequence,
the resulting long exact sequence immediately gives the following implications:
\begin{equation}
\label{implications}
\shF\in\catC\Longrightarrow \shF'\in\catC\quadand \shF',\shF''\in\catC\Longrightarrow\shF\in\catC.
\end{equation}
Furthermore, $\catC$ contains all coherent sheaves with $\dim(\shF)\leq d-1$, according
to \cite{Cutkosky 2014}, Lemma 7.1.
Let $Z_1,\ldots,Z_r\subset Z$ be the irreducible components, endowed with the reduced scheme structure.
By \cite{Cutkosky 2014}, Lemma 10.1 combined with Lemma 9.1 we have
$$
\kappa(\shL)=\max\{\kappa(\shL|Z_1),\ldots,\kappa(\shL|Z_r)\}.
$$
In light of \cite{Cutkosky 2014}, Corollary 9.3 it follows that $\O_{Z_i}\in\catC$, and hence $\O_{Z'}\in\catC$ for every
reduced closed subscheme $Z'\subset Z$. Note that the cited Corollary was formulated for
projective rather than proper schemes, but the proof holds true without changes in the
more general setting.

To proceed, we first suppose that $\shF$ is a torsion-free coherent $\O_{Z_i}$-module, say of rank $r\geq 0$.
Let $\eta\in Z_i$ be the generic point,   choose a bijection $\shF_{\eta}\simeq\kappa(\eta)^{\oplus r}$,
and let $\shF'$ be the resulting intersection   $\shF\cap\O_{Z_i}^{\oplus r}$ inside the 
sheaf $\shM_{Z_i}^{\oplus r}$, where $\shM_{Z_i}$ denotes the quasicoherent sheaf of meromorphic functions.
Then $\shF'$ is coherent and contained in both $\O_{Z_i}^{\oplus r}$ and $\shF$.
The quotient $\shF''=\shF/\shF'$ has dimension $\leq d-1$, thus $\shF''\in\catC$.
Using $\shF'\subset\O_{Z_i}^{\oplus r}$ we infer with \eqref{implications} that $\shF'$ and thus $\shF$ is contained in $\catC$.

Next let $\shF$ be an an $\O_{Z_\red}$-module, and write $\shF_i$ for the restriction $\shF|Z_i$ modulo its torsion subsheaf.
Then we have a short exact sequence
$$
0\lra \shF'\lra \shF\lra \bigoplus_{i=1}^r\shF_i.
$$
The term on the right lies in $\catC$, by the preceding paragraph, and thus also the subsheaf $\shF''=\shF/\shF'$.
The term on the left has $\dim(\shF')\leq d-1$, thus lies in $\catC$, which again by \eqref{implications} gives $\shF\in \catC$.

Finally, let $\shF$ be arbitrary, and $\shI=\Nil(\O_Z)$ be the nilradical, say with $\shI^m=0$.
In  the short exact sequences $0\ra\shI^n\shF/\shI^{n+1}\shF\ra \shF/\shI^{n+1}\shF\ra\shF/\shI^n\shF\ra 0$,
the term on the left is annihilated by $\shI$, whence lies in $\catC$.
Using induction on $n\geq 0$, one sees that the $\shF/\shI^n\shF\in\catC$.  
The case $n=m$ yields $\shF\in\catC$.
\qed

\medskip
It is easy to characterize bigness in dimension one:

\begin{proposition}
\mylabel{big on curves}
Suppose the proper scheme $Z$ is equidimensional, of dimension $d=1$. Then the invertible sheaf $\shL$ is big if and
only if $(\shL\cdot Z')>0$ for some irreducible component $Z'\subset Z$. In particular, this holds if $\deg(\shL)>0$.
\end{proposition}

\proof
We have   
$h^0(\shL^{\otimes n})\geq \chi(\shL^{\otimes n}) = \deg(\shL) n + \chi(\O_Z)$,
where the degree is by definition $\deg(\shL)=\chi(\shL)-\chi(\O_Z)$.
If $\shL$ is not big, Lemma \ref{devissage} implies $\deg(\shL)\leq 0$.
Now suppose that $(\shL\cdot Z')>0$ for some irreducible component $Z'\subset Z$.
By the above, the restriction $\shL|Z'$ is big.
According to \cite{Cutkosky 2014}, Lemma 9.1 combined with Lemma 10.1, the invertible sheaf $\shL$ must be big.
Conversely, assume  $\shL$ is big. Then there is some irreducible component $Z'\subset Z$
such that $\shL|Z'$ is big. By \cite{Cutkosky 2014}, Lemma 7.1
we have $h^0(\shL^{\otimes n}|Z')\geq \alpha n$ for some real constant $\alpha>0$,
and it follows $(\shL\cdot Z')>0$.
\qed

\medskip
We now come to our converse for Proposition \ref{contractible}:

\begin{theorem}
\mylabel{big}
Suppose that $X$ is normal, and let $E\subset X$ be an effective Weil divisor that is contractible
to points, with irreducible components $E_1,\ldots,E_r\subset E$. Suppose $D=\sum m_iE_i$
is a non-zero effective Cartier divisor supported on $E$, and let $\shL=\O_X(-D)$.
Then the restriction $\shL|D$ is big.
\end{theorem}

\proof
It suffices to treat the case that $E$ is connected.
Let $r:X\ra Y$ be the contraction, and $y=r(E)$ be the resulting closed point.
Write  $k=\kappa(y)$ for the residue field, choose a separable closure $k^\sep$ and consider
the resulting geometric point $\bar{y}:\Spec(k^\sep)\ra Y$ and the ensuing  strictly local ring $\O_{Y,\bar{y}}$.
We now replace the scheme  $Y$ by the spectrum of $\O_{Y,\bar{y}}$, and $X$ by the fiber product $X\times_Y\Spec(\O_{Y,\bar{y}})$.
This brings us into the situation that the scheme $Y$ is the spectrum of a strictly local excellent ring $R$,
and $r:X\ra Y$ is a proper morphism with $R=H^0(X,\O_X)$ that is an open embedding on $X\smallsetminus E$
and maps $E$ to the closed point $y\in Y$. Since the formal fibers of the excellent scheme $\Spec(R)$ are geometrically regular,
me may base-change to the formal completion and  assume that the local noetherian ring $R$ is complete.   

Consider the short  exact sequences $0\ra\O_D(-nD)\ra\O_{(n+1)D}\ra\O_{nD}\ra 0$, for each integer $n\geq 0$.
The term on the left is $\shL^{\otimes n}|D$, and we get a short exact sequence
\begin{equation}
\label{infinitesimal}
0\lra H^0(D,\shL^{\otimes n}|D)\lra H^0(X,\O_{(n+1)D}) \lra H^0(X,\O_{nD}).
\end{equation}
The schematic images for the morphisms $nD\ra\Spec(R)$ are of the form $\Spec(R/\ideala_n)$,
for some  inverse system of 
local Artin rings   $R/\ideala_n$. It yields a   descending chain
$\ideala_1\supset\ideala_2\supset\ldots$ of $\maxid_R$-primary ideals given by 
$$
\ideala_n=\{g\in R\mid g\O_X\subset\O_X(-nD)\}.
$$
From this description we see that these ideals 
form a \emph{graded family of ideals} in the sense of \cite{Cutkosky 2014}, that is, $\ideala_m\cdot\ideala_n\subset\ideala_{m+n}$ for all $m,n\geq 0$.
In other words, the subset $\bigoplus \ideala_nT^n\subset R[T]$ is a subring, which one may call the \emph{Rees ring} for
the graded family of ideals.
Since the complete local ring $R$ is reduced, the limit
$$
\alpha=\mult(\ideala_\bullet)=\vol(\ideala_\bullet)=\lim_{n\ra\infty}\frac{\length(R/\ideala_n)}{n^d}
$$
exists as a real number by \cite{Cutkosky 2014}, Theorem 4.7. Here $d=\dim(R)$, and the
number $\alpha$  is called
the \emph{multiplicity} or  \emph{volume}  of the graded family of ideals.
One should think of it as a generalization of the classical \emph{Hilbert--Samuel multiplicities}
$e(\idealb,R)$, which is defined in terms of the graded family of ideal powers $\idealb_n=\idealb^n$.

We now compute this number in two ways.
For the first computation, we describe the ideals $\ideala_n$ in terms of valuations:
Let $x_i\in E_i$ be the generic points.
Since $X$ is normal, the local rings $\O_{X,x_i}$ are discrete valuation rings.
Let $v_i:F^\times\ra\ZZ$ be the corresponding normalized valuation on the field of fractions $F=\Frac(R)=\kappa(\eta)$,
where $\eta\in X$ is the generic point. Then 
$$
\ideala_n=\{g\in R\mid v_1(g)\geq nm_1,\ldots, v_r(g)\geq nm_r \}.
$$
This reveals that the ideals $\ideala_n$ are integrally closed: Indeed,
the codimension one points $x_1,\ldots,x_r\in X$ admit a common affine neighborhood $U=\Spec(A)$,
according to \cite{Gross 2012}, Theorem 1.5.
Write  $\primid_1,\ldots,\primid_r\subset A$ for the corresponding prime ideals of height one.
Then the localization $A'=S^{-1}A$ is a semilocal Dedekind domain, for the multiplicative system 
$S=A\smallsetminus(\primid_1\cup\ldots\cup\primid_r)$. We see
$$
\ideala_n = R\cap (\primid_1^{nm_1}A'\cap\ldots\cap\primid_r^{nm_r}A'),
$$
and this is integrally closed according to \cite{Swanson; Huneke 2006} Proposition 6.8.1 together with  Remark 1.1.3 (8).
Setting $\idealb=\ideala_1$,  we moreover have $\idealb^n\subset \ideala_n$,
and infer that the ideal $\ideala_n$ is the integral closure of the ideal $\idealb^n$.

According to the Brian\c{c}on--Skoda Theorem
in Huneke's form \cite{Huneke 1992}, Theorem 4.13, there is an integer $l\geq 0$ so that
$\ideala_n\subset \idealb^{n-l}$
for all $n\geq l$. Note that this is already a consequence from Izumi's Theorem as given
by H\"ubl and Swanson \cite{Huebl; Swanson 2001}, Theorem 1.2.
It follows that $\length(R/\ideala_n)\geq \length(R/\idealb^{n-l})$.
Passing to the limit, we obtain
$$
\alpha = \lim_{n\ra\infty}\frac{\length(R/\ideala_n)}{n^d} \geq \lim_{n\ra\infty}\left(\frac{\length(R/\idealb^{n-l})}{(n-l)^d} \cdot \frac{(n-l)^d}{n^d}\right).
$$
Indeed, both factors in the sequence on the right converge. The second factor converges to one,
whereas the first factor tends to the Hilbert--Samuel multiplicity $e(\idealb,R)$. 
But such Hilbert--Samuel multiplicities are always integers $e\geq 1$,
according to \cite{AC 8-9}, Chapter VIII, \S4, No.\ 3.
The upshot is that $\alpha\geq 1$.

Seeking a contradiction, we now assume that the restriction $\shL|D$ is not big,
and compute the number $\alpha$ in another way. Recall that  $\dim(X)=\dim(R)=d$, such that $\dim(E)=d-1$.
According to  Lemma \ref{big}, we have $h^0(\shL^{\otimes n}|D)\leq \beta n^{d-2}$ for some
real constant $\beta>0$. By definition of the ideals $\ideala_n$, there are commutative diagrams
$$
\begin{CD}
0	@>>>	\ideala_{n+1}	@>>>	R	@>>> 	H^0(X,\O_{(n+1)D})\\
@.		@VVV			@VVV		@VVV\\
0	@>>>	\ideala_n	@>>>	R	@>>>	H^0(X,\O_{nD})	
\end{CD}
$$
with exact rows. Combining with the exact sequence \eqref{infinitesimal}, we see that
the kernels for the surjection $R/\ideala_{n+1}\ra R/\ideala_n$ are
vector subspaces $V_n\subset H^0(D,\shL^{\otimes n}|D)$.  
Inductively, we infer that 
$$
\length(R/\ideala_n)\leq   \sum_{i=0}^{n-1}\dim(V_i)\leq \sum_{i=0}^{n-1} h^0(\shL^{\otimes i}|E) \leq \beta\sum_{i=0}^{n-1}i^{d-2} \leq \gamma n^{d-1}
$$
for some real constant $\gamma\geq 0$. This in turn gives
$$
\alpha = \lim_{n\ra\infty}\frac{\length(R/\ideala_n)}{n^d} \leq \lim_{n\ra\infty} \frac{\gamma n^{d-1}}{n^d} =0,
$$
contradiction.
\qed

\medskip
If $X$ is a regular 2-dimensional scheme, with a curve $E\subset X$ that is contractible to a point, the intersection matrix
$\Phi=(E_i\cdot E_j)$ is negative-definite, according to Mumford \cite{Mumford 1961} in the complex case,
Artin \cite{Artin 1962} for algebraic surfaces, and Deligne \cite{SGA 7b}, Expos\'e X, Corollary 1.8 in the arithmetic situation.
So for every non-zero effective Cartier divisor $D=\sum m_iE_i$, we have $D^2<0$, and thus the restriction $\shL|D$ 
of the invertible sheaf $\shL=\O_X(-D)$ is big, according to Proposition \ref{big on curves}.
From this point of view, the preceding result can be seen as a generalization  from dimension $d=2$   to higher dimensions.

Now back to our general setting $f:X\ra S=\Spec(R)$.
Let $E\subset X$ be a closed subset that is contractible to points.
If the   proper algebraic space $Y$ resulting from the contraction $X\ra Y$ admits an ample invertible sheaf,
that is, comes from a projective scheme, we 
say that a closed subset $E\subset X$ is \emph{projectively contractible to points}.
This is a rather delicate condition that cannot be determined numerically in general.

The following is a variant of \cite{Schroeer 2017}, Theorem 10.2. 
The new feature is that we have a  \emph{ground ring} rather than a ground field,
and that the contraction is \emph{projective}. This extension will be essential for the
application in the next section.

\begin{proposition}
\mylabel{modification}
Let $E\subset X$ be an effective Cartier divisor that is  contained
in some closed fiber $X_\sigma=f^{-1}(\sigma)$.
Furthermore, suppose that the structure morphism $f:X\ra S$ is projective.
Then there is an effective Cartier divisor $Z\subset E$ with the following property:
On the  blowing-up $g:X'\ra X$ with center $Z\subset X$, the strict transform
$E'\subset X'$ of $E\subset X$ becomes projectively contractible   points.
Moreover, we could choose $Z$ disjoint from any given finite subset $\{x_1,\ldots,x_m\}\subset E$.
\end{proposition}

\proof
Choose a very ample invertible sheaf $\shL$ on the projective scheme $X$ so that
there is a non-zero global section $s_0\in H^0(X,\shL)$ that does not vanish at any of the finitely many points in $\Ass(\O_E)\cup\{x_1,\ldots,x_m\}$.
Then the   map $s_0:\O_E\ra\shL|E$ is injective, and bijective at the points $x_1,\ldots,x_m\in E$.
The section $s_0$ defines an effective Cartier divisor $D\subset X$, and the intersection $Z=D\cap E$ remains Cartier in $E$.
Replacing $\shL$ and $s_0$ by suitable tensor powers, we may assume
that $\shL(-E)$ is   very ample.
Such a closed subscheme $Z\subset X$ is the desired center:

The exceptional divisor for the blowing-up $g:X'\ra X$ is the effective Cartier divisor
$g^{-1}(Z)$. Since $X'$ is integral and $g:X'\ra X$ is dominant, the preimages $g^{-1}(D)$ and $g^{-1}(E)$ are Cartier as well.
Write $ D',E'\subset X'$ for the strict transforms of $D$ and $E$, respectively.
Since the center $Z$ is Cartier on $E$, the universal property of blowing-ups gives
an $X$-morphism $E\ra X'$ whose schematic image is the strict transform $E'$.
In the same way, we have an $X$-morphism $D\ra X'$ with image $D'$. Indeed, for  each point $z\in Z$, let $f_1,f_2\in \O_{X,z}$
be generators for the respective stalks of the ideal sheaves $\O_X(-E),\O_X(-D)\subset\O_X$. By assumption,
they form  a regular sequence. According to \cite{Matsumura 1980}, Theorem 27 on page 98, they remain a regular sequence
in the opposite order, which implies that  the subscheme $Z$ is indeed Cartier in both $E$ and $D$.

According to \cite{Perling; Schroeer 2017}, Lemma 4.4 the strict transforms $D',E'\subset X'$ are Cartier, with  
$$
g^{-1}(D) = D'+g^{-1}(Z)\quadand g^{-1}(E)= E'+g^{-1}(Z)
$$
as Cartier divisors on $X'$. In particular, we have  $\O_{E'}(-E')=\O_E(Z-E)$ with respect to the identification $E'=E$. 
The latter sheaf  is ample on $E$, because the   sheaf $\shL(-E)=\O_X(D-E)$ is   ample on $X$.
By Theorem \ref{contractible}, the Cartier divisor $E'\subset X'$ is contractible to points.
Let  $r:X'\ra \tX$ be the resulting contraction, where $\tX$ is a proper algebraic space.

It remains to construct an ample invertible sheaf on $\tX$.
By the very definition of $X'=\Proj(\bigoplus_{i\geq 0}\shI^i)$  as a relative homogeneous spectrum,
where $\shI\subset\O_X$ is the ideal sheaf for the center $Z\subset X$, we have an invertible sheaf
$\O_{X'}(1)=\O_{X'}(-g^{-1}(Z))$ that  is  relatively ample for the blowing-up $g:X'\ra X$.
Consider the invertible sheaf 
$$
\shL'=g^*(\shL)(1)=\O_{X'}(g^{-1}(D) - g^{-1}(Z)) =\O_{X'}(D').
$$
We claim that $D'$ is disjoint from   $E'$. This is well-known 
(\cite{Hartshorne 1977}, Chapter II, Exercise 7.12), but to fix ideas we provide an argument:
Since $X'\smallsetminus g^{-1}(Z)=X\smallsetminus Z$, it suffices to check that $D'\cap g^{-1}(Z)$ and $E'\cap g^{-1}(Z)$
are disjoint.
The inclusion $Z\subset D$ gives a commutative diagram
$$
\begin{CD}
0	@>>>	\O_X(-D)	@>>>	\O_X	@>>>	\O_D	@>>>	0\\
@.		@VVV			@VV\id V		@VVV\\
0	@>>>	\shI		@>>>	\O_X	@>>>	\O_Z	@>>> 	0,\\
\end{CD}
$$
and the Snake Lemma yields  $0\ra \O_X(-D)\ra\shI\ra\O_D(-Z)\ra 0$. The term on the right
is the ideal sheaf for the Cartier divisor $Z\subset D$, which 
coincides with $\O_D(-E)$. Restricting to $Z$ results in  the short exact sequence
$$
0\lra \O_Z(-D)\lra \shI/\shI^2\lra \O_Z(-E)\lra 0.
$$
Applying this reasoning to the inclusion $Z\subset E$, we infer that the above sequence splits,
and obtain a direct sum decomposition $\shI/\shI^2=\O_Z(-D)\oplus\O_Z(-E)$.
Following Grothendieck's Convention, we regard sections $\Sigma=\sigma(Z)$ for the $\PP^1$-bundle
$$
\PP(\shI/\shI^2)=\Proj\Sym(\shI/\shI^2)=g^{-1}(Z)\lra Z
$$
as invertible quotients $\varphi:\shI/\shI^2\ra\shN$, via the identification $\Sigma =\Proj\Sym(\shN)$.
In the direct sum decomposition $\shI/\shI^2=\O_Z(-D)\oplus\O_Z(-E)$, the first projection
corresponds to the   section  $E'\cap g^{-1}(Z)$, whereas the second projection 
comes from $D'\cap g^{-1}(Z)$. It follows that the two sections are indeed disjoint.
Hence $\shL'$ is trivial in some open neighborhood of $E'$, and consequently $\shL'=r^*(\tshL)$
for some invertible sheaf $\tshL$ on the algebraic space $\tX$.
 
Next, we verify that  $\tshL$ is globally generated.
Since the center $Z\subset X$ is locally of complete intersection, 
we may apply \cite{SGA 6}, Expos\'e VII, Lemma 3.5 together with the Projection Formula
and obtain an identification   $f_*(\shL')=\shI\shL$. 
Since also $r_*(\O_{X'})=\O_{\tX}$, we arrive at the identifications
\begin{equation}
\label{global sections}
H^0(\tX,\tshL)=H^0(X',\shL')=\{ \text{global sections $s$ of $\shL$ with $s_Z=0$}\}.
\end{equation}
The first identification reveals that the base-locus for the invertible sheaf $\tshL$ must be  contained in the image $r(D')$.
The   exact sequence $0\ra\shL(-E)\ra\shL\ra\shL_E\ra 0$ on the original scheme $X$
yields a long exact sequence
\begin{equation}
\label{exact}
0\lra H^0(X,\shL(-E))\stackrel{t_0}{\lra} H^0(X,\shL)\lra H^0(E,\shL_E).
\end{equation}
where $t_0\in H^0(X,\O_X(E))$ is the canonical section defining the inclusion $E\subset X$.
Now recall that the sheaf $\shL(-E)$ very ample and that the Cartier divisor $D\subset X$
is defined by a global section $s_0\in H^0(X,\shL)$.
For each point $x\in D\smallsetminus Z=D\smallsetminus E$ we may choose a global section $t\in H^0(X,\shL(-E))$
with $t(x)\neq 0$. In turn, the global section $s_1= tt_0$ of $\shL$ also has $s_1(x)=t(x)t_0(x)\neq 0$,
thus it defines an effective Cartier divisor $D_1\subset X$ with $x\not\in D_1$ and $Z\subset D_1$.
Under the identification \eqref{global sections}, the resulting section $s_1$ of $\tshL$ does not vanish at $r(x')\in\tX$, where $g(x')=x$.

Now suppose we have a point $z\in Z$. The corresponding point $z'\in g^{-1}(Z)\cap D'$   on $X'$
is an invertible quotient of   $\shI/\shI^2\otimes\kappa(z)=\shI\otimes\kappa(z)$, whence defines a tangent vector at $z\in X$
not contained in $Z$, that is,
a closed subscheme $T\subset\Spec(\O_{X,z})$ of length two,
satisfying $s_0|T=0$ and $t_0|T\neq 0$.
Since $\shL(-E)$ is very ample, we may choose a global section $t$ with $t(z)\neq 0$.
As above, the global section $s_1= tt_0$ of $\shL$ vanishes on $Z$ but not on $T$, thus defines a Cartier divisor  $D_1'\subset X'$ that does not 
contain the point $z'$. Under the identification \eqref{global sections},
the global section $s_1$ of $\tshL$ does not vanish at $r(z')\in\tX$.
Summing up, we have shown that the sheaf $\tshL$ is globally generated.

The last step is to check that the globally generated invertible sheaf $\tshL$ is ample.
Let $\tC\subset \tX$ be an integral curve contained in the fiber $\tX_\sigma$ for the structure morphism
$\tX\ra S$.
We merely have to show that $(\tshL\cdot \tC)>0$.
Let $C'\subset X'$  be the strict transform, such that $C'\not\subset E'$ and $(\tshL\cdot \tC)=(\shL'\cdot C')$.
If $g(C')\subset X$ is a point, that is, $C'$ is a fiber of the $\PP^1$-bundle $g^{-1}(Z)=\PP(\shI/\shI^2)$,
we have $(\shL'\cdot C')=(\O_{X'}(1)\cdot C')\geq 1$.

Now suppose that the image $C=g(C')$ is a  curve rather then a point.
If $C$ is   not contained in the center $Z\subset X$, fix a closed point $x\in C\smallsetminus Z=C\smallsetminus E$
and a global section $t$ for the very ample sheaf  $\shL(-E)$ that vanishes at $x\in C$ but not at the generic point $\eta\in C$.
The resulting global section $s_1=tt_0$ of $\shL$, via the exact sequence \eqref{exact},  
vanishes along $\{x\}\cup Z$ but not at $\eta\in C$. From this we   infer $(\tshL\cdot \tC)>0$.

It remains to treat the case that $ C\subset Z$. Let $\nu:B\ra C'$ be the normalization map,
and form  the fiber product
$$
P=B\times_XX' = B\times_C\PP(\shI/\shI^2\otimes\O_C)  = \PP(\shE),
$$
where the locally free sheaf $\shE=\shL_1\oplus\shL_2$ on $B$ is the sum of the two  invertible sheafs $\shL_1=q^*\O_C(-D)$ and $\shL_2=q^*\O_C(-E)$.
Here  $q:B\ra C$ is the composition of the normalization map $\nu:B\ra C'$ and the induced map $g:C'\ra C$.
The two projections $\pr_i:\shE\ra\shL_i$ correspond to disjoint sections $\Sigma _i\subset \PP(\shE)$,
via $\Sigma _i=\Proj\Sym(\shL_i)$. By functoriality of the construction, these  $\Sigma _1$ and $\Sigma _2$ are the preimages of the sections
$E'\cap g^{-1}(Z)$ and $D'\cap g^{-1}(Z)$ for $g^{-1}(Z)=\PP(\shI/\shI^2)\ra Z$ with respect to  the canonical morphism $B\times_XX'\ra X'$.
In particular, $\O_{P}(\Sigma_2)$ is the pullback of $\shL'$.

The scheme $P=\PP(\shE)$ is a ruled surface over the proper regular curve $B$, 
so its Picard group modulo numerical equivalence takes the form $N(P)=\ZZ^2$.
The pseudoeffective cone 
$$
\overline{\operatorname{NE}}(P)\subset N(P)\otimes_\ZZ\RR=\RR^2
$$
must be  generated by two extremal rays. Each fiber $F\subset P$ for the ruling has selfintersection number  $F^2=0$.
According to Proposition \ref{pullback}, the section $\Sigma_1\subset P$ is contractible, so its selfintersection is $(\Sigma_1)^2<0$.
In light of \cite{Kollar 1995}, Lemma 4.12, it follows that the numerical classes of $F$ and $\Sigma_1$ 
are the two extremal rays for  $\overline{\operatorname{NE}}(P)$.
This in turn implies  $(\Sigma_2)^2>0$. In particular, $\O_P(\Sigma_2)$ is ample on $\Sigma_2$.
According to Fujita's result (\cite{Fujita 1983}, see also  \cite{Ein 2000}), the invertible sheaf $\O_P(\Sigma_2)$ must be semiample.
It follows that for some $n\geq 1$ the semiample sheaf $\O_P(n\Sigma_2)$ 
is the preimage of some ample sheaf on  $\tP$, where $P\ra \tP$ is the contraction of $\Sigma_1$.
Consequently $(\Sigma_2\cdot \Sigma)>0$ for every integral curve $\Sigma\neq \Sigma_1$.
In particular, this holds for the section $\Sigma\subset P$ arising from the diagonal map 
$B\ra  B\times_X X'=P$. By construction, the projection $P\ra X'$ induces a surjection $B=\Sigma\ra C'$.
Since $\O_P(\Sigma_2)$ is the preimage of $\shL'=\O_{X'}(D')$, we must have $(\shL'\cdot C')>0$.
\qed

%===========================================================
\section{Totally separably closed schemes}
\mylabel{TSC}

Recall that an integral   scheme $X$ with generic point $\eta\in X$  
is   \emph{totally separably closed} if it is normal and the function field  $F=\O_{X,\eta}=\kappa(\eta)$ is separably closed.
A space or a scheme is called \emph{local} if it contains exactly one closed point.
The main result of this paper is:

\begin{theorem}
\mylabel{local}
Let $X$ be an integral  separated scheme that is totally separably closed, and 
$u,v\in X$ be two points. Then the intersection $\Spec(\O_{X,u})\cap\Spec(\O_{X,v})$
inside $X$ is   local.
\end{theorem}

\proof
The intersection can be regarded as the underlying set of the schematic fiber product 
$$
P=\Spec(\O_{X,u})\times_X\Spec(\O_{X,v}).
$$
Its image contains the generic point $\eta\in X$,
in particular $P$ is non-empty. Furthermore, the scheme $P$ is affine, because $X$ is separated.
Seeking a contradiction, we assume that the intersection is not local.
Hence there are two closed points $\alpha\neq\beta$ inside $P$.
Let $A,B\subset X$ be their closures in $X$. Both contain $u$ and $v$. In fact, the points $u,v\in A\cap B$ are
generic points in the intersection.
If the two points $u,v\in X$ admit a common affine open neighborhood,
we immediately get a contradiction from \cite{Artin 1971}, Corollary 1.8.
The idea of this proof is to construct, starting form $X$, another integral separated scheme $\tilde{X}$ that is totally separably closed
and additionally  enjoys the AF property, containing two points $\tu\neq\tv$ closely related to the original points $u\neq v$.
Now the intersection $\tilde{P}=\Spec(\O_{\tX,\tu})\times_{\tX}\Spec(\O_{\tX,\tv})$ is indeed local, and this will
finally produced the desired contradiction.

\medskip
{\bf Step 1:}
\emph{We reduce to the case that $X$ is the total separable closure of some proper $\ZZ$-scheme $X_0$.}
First of all, we may assume that our scheme $X$ is quasicompact, simply by choosing affine open neighborhoods
$U,V\subset X$ of $u,v\in X$ and replacing $X$ by their union.
Next, we write $X=\invlim X_\lambda$ as a filtered inverse limit of
schemes $X_\lambda$, $\lambda\in L$ that are separated and of finite type over the ring $R=\ZZ$, with affine transition maps $X_\mu\ra X_\lambda$,
$\lambda\leq\mu$. This is possible according to \cite{Thomason; Trobaugh 1990}, Appendix C, Proposition 7.
Replacing $X_\lambda$ by the schematic images of the projection $X\ra X_\lambda$, we may
assume that the $X_\lambda$ are integral, and that the transition maps $X_\mu\ra X_\lambda$ and the projections $X\ra X_\lambda$
are dominant. Let $\eta_\lambda\in X_\lambda$ be the generic points, such that $\eta=(\eta_\lambda)_{\lambda\in L}$.
The function fields $F_\lambda=\O_{X_\lambda,\eta_\lambda}=\kappa(\eta_\lambda)$ form a filtered direct system
of subfields inside $F=\O_{X,\eta}=\kappa(\eta)$, with $F=\bigcup_{\lambda\in L} F_\lambda$.
Let $F^\sep_\lambda$ be the relative separable algebraic closure of $F_\lambda\subset F$, and $\TSC(X_\lambda)$ the  be the resulting
integral closure of $X_\lambda$ with respect to the field extension $F_\lambda\subset F^\sep_\lambda$. Then the filtered inverse system
$X_\lambda$ induces a filtered inverse system $Y_\lambda=\TSC(X_\lambda)$. The morphism $X\ra X_\lambda$ induces compatible morphisms
$X\ra Y_\lambda$, giving and identification $X=\invlim Y_\lambda$.

Now suppose that the theorem is valid for all $Y_\lambda=\TSC(X_\lambda)$.
Let $u_\lambda,v_\lambda\in Y_\lambda$ be the images of $u,v\in X$.
Then the schemes
$$
P_\lambda=\Spec( \O_{Y_\lambda,u})\times_{Y_\lambda} \Spec(\O_{Y_\lambda,v})
$$
are local henselian. According to \cite{Artin 1971}, Lemma 2.6
the inverse limit $P=\invlim P_\lambda$ 
is  local henselian, contradiction.

This   reduces us to the case that $X$ is the total separable closure of some integral 
scheme $X_0$ that is separated and of finite type over the ring $R=\ZZ$.
In light of Nagata's Compactification Theorem in the relative version obtained by L\"utkebohmert   \cite{Luetkebohmert 1993},
we may additionally assume that the structure morphism $X_0\ra\Spec(\ZZ)$ is proper. This concludes Step 1.

\medskip
{\bf Step 2:}
\emph{We may assume that the point $u\in X$ lies in a closed fiber $X\otimes\FF_p$.}
If both points $u,v\in X$ lie in the generic fiber $X\otimes\QQ$, we could replace $X_0$ and $X$
by their generic fibers. Now $X_0$ is proper over the field $k=\QQ$, and we  immediately get a contradiction to \cite{Schroeer 2017}, Theorem 12.1.
So we may assume without restriction that $u\in U$ lies in a closed fiber $X\otimes\FF_p$, for some prime $p>0$.

\medskip
{\bf Step 3:}
\emph{Here we write $X$ as a filtered inverse system $X_\lambda$, $\lambda\in L$ of proper $\ZZ$-schemes so that the geometry of 
$A,B\subset X$  is captured by their images  $A_\lambda,B_\lambda\subset X_\lambda$.}
Let $F_0\subset F$ be the inclusion of function fields coming from the canonical morphism $X\ra X_0$.
Changing the notation from Step 1, we now write  $F_\lambda\subset F$, $\lambda\in L$ for   the filtered direct system
of subfields with $[F_\lambda:F_0]<\infty$, and let $X_\lambda\ra X$ be the normalization of $X_0$ with
respect to the field extension $F_0\subset F_\lambda$. This gives a filtered inverse system $X_\lambda$ of finite $X_0$-schemes,
with $X=\invlim X_\lambda$, where the transition maps $X_\mu\ra X_\lambda$, $\lambda\leq \mu$ are finite.
Note that all structure morphisms $X_\lambda\ra\Spec(\ZZ)$ are proper,   that $X\ra\Spec(\ZZ)$ is 
separated and universally closed, and that the projections $X\ra X_\lambda$ are universally closed. Write
\begin{equation}
\label{induced points}
u_\lambda, v_\lambda, \alpha_\lambda,\beta_\lambda, \eta_\lambda\in X_\lambda
\end{equation}
for the respective images of the   points $u,v,\alpha,\beta,\eta\in X$. Let $A_\lambda,B_\lambda\subset X_\lambda$
be the closures of $\alpha_\lambda,\beta_\lambda\in X_\lambda$, which are also the images of $A,B\subset X$.
In turn, we have $u_\lambda,v_\lambda\in A_\lambda\cap B_\lambda$.
Replacing $L$ by some cofinal subset, we may assume that the points in \eqref{induced points} are pairwise different.

Since inverse limits commute with inverse limits, we have $A\cap B=\invlim (A_\lambda\cap B_\lambda)$.
For the local rings this means $\O_{A\cap B,u} = \dirlim \O_{\O_{A_\lambda\cap B_\lambda, u_\lambda}}$.
Setting $C=\Spec(\O_{A\cap B,u})$ and $C_\lambda = \Spec(\O_{\O_{A_\lambda\cap B_\lambda, u_\lambda}})$, we get 
$C\smallsetminus\{u\}=\invlim(C_\lambda\smallsetminus\{u_\lambda\})$. But the left-hand side is empty,
because $u\in A\cap B$ is a generic point. Thus   $1=0$ already holds as global sections on some $C_\lambda$.
By symmetry, the same applies for the point $v\in A\cap B$. Replacing $L$ by some cofinal subset,
we thus may assume that $u_\lambda,v_\lambda\in A_\lambda\cap B_\lambda$ are generic points.

\medskip
{\bf Step 4:}
\emph{Reduction to the case that the connected components of $A_0\cap B_0$ are irreducible.}
The proper $\ZZ$-scheme $A_0\cap B_0$ has only finitely many irreducible components.
Fix an irreducible component $C_0\subset A_0\cap B_0$, let $C'_0\subset A_0\cap B_0$ be the union
of the other irreducible components, and consider the normalized blowing-up $Y_0\ra X_0$ with center $Z=C\cap C'$.
This morphism is an isomorphism over some open neighborhood of the set of generic points in $A_0\cap B_0$,
and the strict transforms of $C_0$ and $C'_0$ become disjoint on $Y_0$. 
By   induction on  the number of irreducible components in $A_0\cap B_0$, we find a normalized blowing-up $Y_0'\ra X_0$
that is an isomorphism over the set of generic points of $A_0\cap B_0$  so that the strict transform $A_0', B'_0$
have the property that the connected components of  $A'_0\cap B'_0$ are irreducible.
 Replacing $X_0$ by $Y'_0$ 
and $X$ by $\TSC(Y'_0)$ we may assume that the connected components of $A_0\cap B_0$ are irreducible.
In particular, 
\begin{equation}
\label{intersection empty}
\overline{\{u_0\}}\cap \overline{\{v_0\}}=\emptyset.
\end{equation}
Note that his property will later produce the desired contradiction.

\medskip
{\bf Step 5:}
\emph{Construction of two auxiliary filtered inverse systems $X'_\lambda$ and $X''_\lambda$ consisting
of projective schemes.}
Choose an affine open neighborhood $V_0\subset X_0$ of $v_0$ not containing $u_0$.
By Chow's Lemma, there is a blowing-up $g_0:X'_0\ra X_0$ with center $Z_0$ contained in
$X_0\smallsetminus V_0$, such that $X'_0$ is projective.
Set $Z'_0=g_0^{-1}(\overline{\{u_0\}})$. 
Replacing $X'_0$ by some further normalized blowing-up of $X'_0$ with center $Z_0'$, we may also assume
that the   the   closed set  $g_0^{-1}(\overline{\{u_0\}})$ is the support of an 
effective Cartier divisor  $E'_0\subset X'_0$ contained in the closed fiber $X'_0\otimes\FF_p$. 

Let $X'_\lambda\ra X'_0$ be its normalization of $X'_0$ with respect to the finite field extension $F_0\subset F_\lambda$. 
This gives a filtered direct system with
finite transition maps, and we obtain a totally separably closed scheme $X'=\invlim X'_\lambda$.
The projection $X'\ra X'_0$ is affine  and the scheme $X'_0$ satisfies the 
the AF property, so the same holds for $X'$.
The ensuing  projective birational morphisms $X'_\lambda\ra X_\lambda$ induces
a  birational morphism $X'\ra X$ between integral schemes. 
The   $X'_\lambda\ra X_\lambda$ are isomorphisms over the open subsets $V_\lambda=V_0\times_{X_0}X_\lambda$,
which contains $v_\lambda,\alpha_\lambda,\beta_\lambda$, hence $X'\ra X$
is an isomorphism over $V=V_0\times_{X_0}X$, which contains $v,\alpha,\beta$.
So we may regard the latter also as points  $v',\alpha',\beta'\in X'$.

Let $A',B'\subset X'$ be the strict transforms of $A,B\subset X$, that is, the closures of $\alpha',\beta'\in X'$.
According to \cite{Artin 1971}, Corollary 1.8 together with
\cite{Schroeer 2017}, Theorem 7.6 the intersection $A'\cap B'$ is irreducible. 
Since $\Spec(\O_{A'\cap B',v'})=\Spec(\O_{A\cap B,v})$, the point $v'\in A'\cap B'$ is generic, 
whence this   must be the unique generic point.
Consider the canonical morphism $\varphi:X'\ra X_0$, which is a closed map.
Suppose there would be a point $u'$ in the intersection
\begin{equation}
\label{empty fiber}
 A'\cap B'\cap \varphi^{-1}(u_0)=\overline{\{v'\}} \cap \varphi^{-1}(u_0).
\end{equation}
Then 
$u_0\in \varphi(\overline{\{v'\}}) =  \overline{ \{\varphi(v')\}}=\overline{ \{v_0\}}$,
contradicting \eqref{intersection empty}. Thus the intersection \eqref{empty fiber} is empty.
Passing to a cofinal subset of $L$, we   may assume that already the fiber of $A'_0\cap B'_0$
over $u_0\in A_0\cap B_0$ is empty. On the other hand, the morphisms  $A'_0\ra A_0$ and $B'_0\ra B_0$ are
proper and dominant, so there are points $r'_0\in A'_0$ and $s'_0\in B'_0$ lying in $\varphi^{-1}(u_0)$.

Now recall that the fiber $g_0^{-1}(\overline{\{u_0\}})$ is  a Cartier divisor $E'_0\subset X'_0$
lying inside the closed fiber $X'\otimes\FF_p$.
Write $n\geq 2$ for the common dimension of the schemes $X,X_\lambda, X',X'_\lambda$.
Since $A_0,B_0\subsetneq X_0$, the irreducible schemes $A'_0,B'_0$ have dimension $\leq n-1$,
and  the intersections $A'_0\cap E'_0$ and $B'_0\cap E'_0$ have dimension $\leq n-2$.
On the other hand, each irreducible component of $E'_0$ has dimension $n-1$,
by Krull's Principal Theorem.
According to Proposition \ref{modification}, there is an effective Cartier divisor $Z_0'\subset E'_0$ 
not containing $\{r'_0,s'_0\}$ so that on the blowing-up $X''_0\ra X'_0$
with center $Z'_0\subset X'_0$  the strict transform $E''_0$ of $E'_0$
becomes projectively contractible to   points.
As above, we get a filtered inverse system
of projective schemes $X''_\lambda$ with finite transition maps,
and the resulting $X''=\invlim(X''_\lambda)=\TSC(X''_0)$ comes with a birational morphism $X''\ra X'$
of integral schemes. Consider the composite morphism $h: X''\ra X$ and the corresponding
$h_\lambda:X''_\lambda\ra X_\lambda$. By construction, we have $h_0(E''_0)=\overline{\{u_0\}}$.
This finishes Step 5.

\medskip
{\bf Step 6:}
\emph{Construction of the contractions $r_\lambda:X''_\lambda\ra\tX_\lambda$ resulting in  a contradiction.}
Let $E''_\lambda\subset X''_\lambda$ be the preimages of $E''_0\subset X''_0$.
These subschemes are effective Cartier divisors, because $X''_\lambda$ is integral and $X''_\lambda\ra X''_0$
is dominant.
According to Proposition \ref{pullback}, the   $E''_\lambda$ are projectively contractible to points.
Moreover, the resulting contractions $r_\lambda:X''_\lambda\ra\tX_\lambda$ coincide with the Stein factorization of
$X''_\lambda\ra \tX_0$, and   yield yet another filtered inverse system
$\tX_\lambda$ of projective schemes with finite transition maps. We  have a commutative diagram
$$
\begin{CD}
\tX_\lambda 		@<r_\lambda << 	X''_\lambda	@>h_\lambda>>	X_\lambda\\
@VVV					@VVV				@VVt_{0\lambda} V\\
\tX_0			@<<r_0<		X''_0		@>>h_0>		X_0.
\end{CD}
$$
Clearly, $u_\lambda$ is contained in the finite
set $t_{0\lambda}^{-1}(u_0)$, where $t_{0\lambda}:X_\lambda\ra X_0$ denotes the transition map. Using the above commutative diagram, we infer that
the    $h_\lambda^{-1}(u_\lambda)$ is contained in the effective Cartier divisor $E''_\lambda$
whose connected components are mapped to closed points in $\tX_\lambda$.
But the fiber $h_\lambda^{-1}(u_\lambda)$ is connected, by Zariski's Main Theorem.
It follows that the set $h_\lambda^{-1}(u_\lambda)$  and in particular the elements
$r_\lambda,s_\lambda\in h_\lambda^{-1}(u_\lambda)$ are mapped to the same  point $\tu_\lambda\in\tX_\lambda$.

Now consider the resulting filtered inverse system $\tX_\lambda$, $\lambda\in L$ of projective schemes with finite
transition maps. The inverse limit $\tX=\invlim \tX_\lambda$ is another totally separably closed scheme.
Since the scheme $\tX_0$ is projective and the morphisms $\tX\ra\tX_0$ is integral,  $\tX$ enjoys the AF property.
The points $\tu_\lambda\in \tX_\lambda$ are compatible and yield a point $\tu\in \tX$.
By construction,   $X_\lambda\ra\tX_\lambda$ are isomorphisms on an open neighborhood of $v_\lambda,\alpha_\lambda,\beta_\lambda$.
So we may regard the latter as points on $\tX$, denoted by 
$\tv_\lambda,\tilde{\alpha}_\lambda,\tilde{\beta}_\lambda\in\tX_\lambda$. In turn, we get points 
$\tv,\tilde{\alpha},\tilde{\beta}$ on the inverse limit $\tX$.

Since the points $\tu,\tv\in \tX$ admit a common affine open neighborhood,
\cite{Artin 1971}, Corollary 8.1 applies and we infer with \cite{Schroeer 2017}, Theorem 7.6
that the intersection $\tA\cap \tB$ is irreducible. Arguing as above, we see that $\tv\in \tA\cap \tB$
must be the generic point, in particular $\tu\in\overline{\{\tv\}}$.
Since the projection $\tX\ra \tX_0$ is closed, we also have $\tu_0\in\overline{ \{\tv_0\} }$.
The contraction $r_0:X''_0\ra\tX_0$ is   closed as well, hence
$$
r_0(E''_0\cap \overline{\{v''_0\}}) \supset
r_0( r_0^{-1}(\tu_0 ) \cap \overline{\{v''_0\}}) =  
\{\tu_0\} \cap r_0( \overline{ \{v''_0\} }) =
\{\tu_0\} \cap  \overline{ \{\tv_0\} }. 
$$
We see that $ E''_0\cap \overline{ \{v''_0\} }\neq \emptyset$.
Finally, the blowing-up  $h_0:X''_0\ra X_0$ is a closed map with $h_0(E''_0)=\overline{\{u_0\}}$ and $h_0(v''_0)=v_0$, hence the sets
$$
h_0(E''_0\cap \overline{\{v''_0\}}) \subset
h_0(E''_0) \cap h_0(\overline{\{v''_0\}})=
\overline{\{u_0\}}\cap \overline{\{v_0\}}
$$
are non-empty. But this contradicts \eqref{intersection empty}.
\qed

\medskip
Note that  in Theorem \ref{local}  some   assumption about separatedness  is inevitable :
For example, let $k$ be an algebraically closed  field, $R_0$ be the henselization
of $k[x,y]$ at the maximal ideal $\maxid=(x,y)$, and $R$ be its total separable closure.
Then $R$ is a   local integral domain of dimension two that is TSC.
Let $U\subset \Spec(R)$ be the complement
of the closed point. Then  $U$ has dimension one and   contains infinitely many
closed points. Let $R_1$ and $R_2$ be two copies of $R$, and 
$X=\Spec(R_1)\cup\Spec(R_2)$ be the non-separated integral TSC scheme obtained by gluing along $U\subset\Spec(R_i)$.
For the two closed points $u,v\in X$ we have $\Spec(\O_{X,u})\cap\Spec(\O_{X,v})=U$,
which is not local.

In this example, the diagonal $\Delta:X\ra X\times X$ is not affine.
It is conceivable that Theorem \ref{local} holds true under the weaker assumption that the diagonal is merely affine rather than a closed embedding.

%===========================================================
\section{Application to Nisnevich cohomology}
\mylabel{Nisnevich}

Let $X$ be a scheme, and write $(\Et/X)$ for the category of \'etale $X$-schemes.
The \emph{Nisnevich topology} on this category is the Grothendieck topology 
defined by the pretopology whose covering families $(U_i\ra U)_{i\in I}$
are those   where    for each $x\in U$ there is some index $i\in I$ and some $x_i\in U_i$ mapping to $x$,
such that the residue field extension $\kappa(x)\subset\kappa(x_i)$ is trivial 
\cite{Nisnevich 1989}.
We write $X_\Nis$ for the ensuing topos of presheaves on $(\Et/X)$ that satisfy
the sheaf axiom for the Nisnevich topology. We refer to such sheaves as \emph{Nisnevich sheaves}.
Each point $x\in X$ yields a point  $(P_*,P^*,\psi):(\Set)\ra  X_\Nis $ in the sense of topos-theory,
and the corresponding local ring of the structure sheaf with respect to the Nisnevich topology 
is the henselization of the local ring $\O_{X,x}$ with respect to the Zariski topology.
Every abelian Nisnevich sheaf $F$  comes with a spectral sequence
$$
E_2^{pq}=\cH^p(X_\Nis,\underline{H}^q(F)) \Longrightarrow H^{p+q}(X_\Nis,F)
$$
from \v{C}ech cohomology to sheaf cohomology (see for example
\cite{Schroeer 2017}, Appendix B).

\begin{theorem}
\mylabel{cech}
If $X$ is quasicompact and separated, then $\cH^p(X_\Nis,\underline{H}^q(F)) =0$
for all $p\geq 0$, $q\geq 1$ and all abelian Nisnevich sheaves $F$.
In particular, the canonical maps 
$$
\cH^p(X_\Nis,F)\lra H^p(X_\Nis,F)
$$
from \v{C}ech cohomology to sheaf cohomology are bijective for all $p\geq 0$.
\end{theorem}

\proof
The result was already established in  \cite{Schroeer 2017}, Theorem 13.1 for schemes
where the structure morphism $X\ra\Spec(\ZZ)$ factors over the spectrum of a prime field.
In other words,   $X$ is a $k$-scheme for some ground field $k$.
This assumption entered   only via \cite{Schroeer 2017}, Theorem 12.1. But the latter holds true
without the superfluous assumption of a ground field, by Theorem \ref{local}.
\qed

%===========================================================


\begin{thebibliography}{ccccc}

\bibitem{Artin 1962}
M.\ Artin:
Some numerical criteria for contractability of curves on algebraic
surfaces.
Am.\ J.\ Math.\ 84 (1962), 485--496.

\bibitem{Artin 1970}
M.\ Artin:
Algebraization of formal moduli II: Existence of modifications.
Ann.\ Math.\  91 (1970), 88--135.

\bibitem{Artin 1971}
M.\ Artin:
On the joins of Hensel rings.
Advances in Math.\ 7 (1971),  282--296.

\bibitem{Benoist 2013}
O.\ Benoist:
Quasi-projectivity of normal varieties. 
Int.\ Math.\ Res.\ Not.\ IMRN  17 (2013), 3878--3885. 

\bibitem{SGA 6}
P.\ Berthelot, A.\ Grothendieck, L.\ Illusie (eds.):
Th\'eorie des intersections et th\'eor\`eme de Riemann--Roch (SGA 6).
Springer, Berlin, 1971.

\bibitem{AC 8-9}
N.\ Bourbaki:
Alg\`ebre commutative. Chapitre 8--9.
Masson, Paris, 1983.

\bibitem{Cutkosky 2014}
S.\ Cutkosky:
Asymptotic multiplicities of graded families of ideals and linear series. 
Adv.\ Math.\ 264 (2014), 55--113. 

\bibitem{SGA 7b}
P.\ Deligne,   N.\ Katz:
Groupe de monodromie en g\'eom\'etrie alg\'ebrique (SGA 7 II).
Springer, Berlin, 1973.

\bibitem{Ein 2000}
L.\ Ein:
Linear systems with removable base loci.
Comm.\ Algebra 28 (2000),5931--5934. 

\bibitem{Engler and Prestel 2005}
A.\ Engler, A.\ Prestel:
Valued fields. 
Springer, Berlin, 2005. 
 
\bibitem{Fujita 1983}
T.\ Fujita:
Semipositive line bundles. 
J.\ Fac.\ Sci.\ Univ.\ Tokyo   30 (1983), 353--378.

\bibitem{Gross 2012}
P.\ Gross:
The resolution property of algebraic surfaces.  
Compos.\ Math.\ 148 (2012),  209--226.

\bibitem{EGA II}
A.\ Grothendieck:
\'El\'ements de g\'eom\'etrie alg\'ebrique II:
\'Etude globale \'el\'ementaire de quelques classes de morphismes.
Publ.\ Math., Inst.\ Hautes \'Etud.\ Sci.\ 8 (1961).

\bibitem{EGA IIIa}
A.\ Grothendieck:
\'El\'ements de g\'eom\'etrie alg\'ebrique III:
\'Etude cohomologique des faiscaux coh\'erent.
Publ.\ Math., Inst.\ Hautes \'Etud.\ Sci.\ 11 (1961).

\bibitem{Huneke 1992}
C.\ Huneke:
Uniform bounds in Noetherian rings. 
Invent.\ Math.\ 107 (1992), 203--223.

\bibitem{Huneke 2011}
C.\ Huneke:
Absolute integral closure. 
In:  A.\ Corso and C.\ Polini (eds.), Commutative algebra and its connections to geometry, pp.\ 119--135.
Amer. Math. Soc., Providence, RI, 2011. 

\bibitem{Hartshorne 1977}
R.\ Hartshorne:
Algebraic geometry.
Springer, Berlin,  1977.

\bibitem{Huebl; Swanson 2001}
R.\ H\"ubl, I.\ Swanson:
Discrete valuations centered on local domains.
J.\ Pure Appl.\ Algebra 161 (2001), 145--166.

\bibitem{Kollar 1995}
J.\ Koll\'ar:
Rational curves on algebraic varieties.
Springer, Berlin, 1995.

\bibitem{Lazarsfeld 2004}
R.\ Lazarsfeld:
Positivity in algebraic geometry. I.  
Springer, Berlin, 2004.

\bibitem{Luetkebohmert 1993}
W.\ L\"utkebohmert:
On compactification of schemes.
Manuscr.\ Math.\ 80 (1993),  95--111.

\bibitem{Matsumura 1980}
H.\ Matsumura:
Commutative algebra. Second edition. 
Benjamin/Cummings, Reading, Mass., 1980.

\bibitem{Mumford 1961}
D.\ Mumford:
The topology of normal singularities of an algebraic surface and a criterion for simplicity.
Publ.\ Math., Inst.\ Hautes \'Etud.\ Sci.\ 9 (1961), 5--22.

\bibitem{Nisnevich 1989}
Y.\ Nisnevich:
The completely decomposed topology on schemes and associated descent spectral sequences in algebraic K-theory. 
In: 
J.\ Jardine, V.\ Snaith (eds.), Algebraic K-theory: connections with geometry and topology, pp.\  241--342.
Kluwer, Dordrecht, 1989.
 
\bibitem{Olsson 2016}
M.\ Olsson:
Algebraic spaces and stacks.
American Mathematical Society, Providence, RI, 2016.

\bibitem{Perling; Schroeer 2017}
M.\ Perling, S.\ Schr\"oer:
Vector bundles on  proper toric 3-folds and certain other   schemes.
Trans.\ Amer.\ Math.\ Soc.\ 369 (2017), 4787--4815.

\bibitem{Schmidt 1933}
F.\ Schmidt:
K\"orper, \"uber denen jede Gleichung durch Radikale aufl\"osbar ist.  
Sitzungsber.\ Heidelberger Akad.\ Wiss.\ 2  (1933), 37--47.

\bibitem{Schroeer 2017}
S.\ Schr\"oer:
Geometry on totally separably closed schemes.
Algebra Number Theory 11 (2017),  537--582.

\bibitem{Swanson; Huneke 2006}
I.\ Swanson, C.\ Hueke:
Integral Closure of Ideals, Rings, and Modules.
Cambridge University Press, Cambridge, 2006.

\bibitem{Thomason; Trobaugh 1990}
R.\ Thomason, T.\ Trobaugh:
Higher algebraic $K$-theory of schemes and of derived categories. 
In: 
P.\ Cartier et al.\ (eds.),
The Grothendieck Festschrift III, 247--435.
Birkh\"auser, Boston, MA, 1990.

\end{thebibliography}
\end{document}